\documentclass[a4paper,12pt]{article}
\begin{document}

\title{Solutions of nonlinear Sobolev-Burgers PDEs.}
\author{S.V. Ludkovsky}

\date{28 March 2018}
\maketitle

\begin{abstract}
Nonlinear Sobolev-Burgers PDEs are considered. Their solutions are
investigated. A technique of noncommutative line integration is
utilized for their description. A new method of PDEs solution with
the help of Cayley-Dickson algebras is developed in the article.
Moreover, random operator valued measures are studied and applied to
solutions of PDEs.

\footnote{key words and phrases: nonlinear; partial differential
equation; operator; noncommutative integration; complex  \\
Mathematics Subject Classification 2010: 35G20; 35L75; 32W50; 30G35;
60G20; 60H05; 17A05; 17A45
\\ address: S.V. Ludkovsky, Dep. Appl. Mathematics, Moscow State Techn. Univ. MIREA,
\\ av. Vernadsky 78, Moscow 119454, Russia \\ sludkowski@mail.ru }

\end{abstract}

\section{Introduction.}
Studies of nonlinear PDEs compose an extensive part of nonlinear
analysis and PDE (see, for example,
\cite{ablsigb,hormbl,kichenb96,polzayjurb,polzayrb,svalkorplb} and
references therein). There are linear and nonlinear PDEs of
particular types, for example, wave PDEs, heat PDEs, diffusion PDEs,
Schr\"odinger PDEs to each of which more than a thousand articles
and books are devoted. Though about an existence of solutions and
their numerical simulations a lot of works is written, their
integration remains a hard or an unsolved problem for many types of
PDEs. On the other side, random functions methods based on Gaussian
measures appeared to be useful for solutions of diffusion PDEs and
Schr\"odinger PDEs (see, for example,
\cite{gihmskorb,gulcastb,johnlap,melnb,nicolam16} and references
therein).
\par Among nonlinear PDEs the Sobolev type PDEs, which include
the Burgers PDE, have important applications in physics and
hydrodynamics \cite{svalkorplb}. They describe two-dimensional
motions of a stratified rotating liquid, electromagnetic fields in
crystals, internal gravitational waves, non stationary filtration
process of liquid in a fissure porous medium, dissipation process,
cold plasma, two temperature plasma in an external magnetic field,
etc.  (see \cite{svalkorplb} and references therein). There were
works about numerical computation of some PDE solutions and an
existence of solutions. Nevertheless needs arise to integrate such
type and more general PDEs or their systems. In many cases it is
also necessary to analyze properties of solutions. Therefore
analytic approaches apart from that of numerical provide in this
respect many advantages. \par Such problems led to a new main stream
of developing hypercomplex analysis for the PDE theory needs (see
\cite{br} - \cite{frenonlutz2014}, \cite{guhaspr} - \cite{guetze},
\cite{ludfov}-\cite{ludrimut2014} and references therein). It was
begun mainly in the years 1990-th over quaternions and Clifford
algebras. Later on since the years 2010-th it was begun over
octonions and more general Cayley-Dickson algebras. A reason for
such activity is in an enlargement of possibilities: some PDEs which
are not integrable over the complex field $\bf C$ appear to be
integrable over the aforementioned algebras so that a suitable
Clifford or a Cayley-Dickson algebra can be chosen for a given PDE.
\par Previously a new approach of noncommutative integration of nonlinear PDEs
was investigated over hypercomplex numbers
\cite{lucveleq2013,frelud17}. It was applied to PDEs used in
hydrodynamics such as the non-isothermal flow of a non-compressible
Newtonian liquid PDE and the Korteweg-de-Vries PDE. In this paper
also another new mathematical tool is developed and applied to
solutions of PDEs, which consists in random operator valued
measures.
\par Some readers of previous articles asked how the Cayley-Dickson algebras
can be used for solutions of general type PDEs. To this problem
Section 2 is devoted. In this paper such theory is developed further
for other types of PDEs which were not yet treated. Procedures
permitting to write equivalent problems over quaternions, octonions
or Cayley-Dickson algebras instead of PDEs over $\bf R$ or $\bf C$
are described. The corresponding Theorems 1, 2 and Proposition 3 are
proved. The method is illustrated on Sobolev-Burgers PDEs in Section
3. There solutions of PDEs are investigated. For this purpose random
operator valued measures are studied. The noncommutative line
integration is utilized for their description. Solutions of the
Sobolev-Burgers PDEs are investigated (see Theorems 6, 23 and 24 and
Subsection 25 below).
\par  All main results of this paper are obtained
for the first time.  They can be used for further studies and
integrations of linear and nonlinear PDEs.

\section{Solution of PDEs over octonions.}
\par Frequently PDEs are given over the real field or the complex
field. For an application of the noncommutative integration
technique of PDEs it may be necessary at first to present the
corresponding PDEs over octonions or Cayley-Dickson algebras.
Henceforth notations and definitions of the work \cite{lucveleq2013}
are used.

\par {\bf 1. Theorem.} {\it To each
scalar or vector PDE \par $(1.1)$ $P(A_1,...,A_v,u)=g$ over $\bf R$
which can be nonlinear relative to an unknown scalar or vector
function $u=(u_1,...,u_m)$, where $A_1,...,A_v$ are linear PDO with
real coefficients, $g=(g_1,...,g_k)$, $g_s$ and $u_j$ are real
functions for each $s=1,...,k$ and $j=1,...,m$, $1\le k\in \bf N$,
$1\le m\in \bf N$, $P=(P_1,...,P_k)$ and $g$ are given, $P_s$ is a
polynomial (or power series) with real coefficients for each $s$, a
domain $U$ is open
in a canonically closed subset $V$ in ${\bf R}^n$, \\
can be posed a PDE
\par $(1.2)$ $Q(\hat{A}_1,...,\hat{A}_v,\hat{u})=\hat{g}$ over the
Cayley-Dickson algebra ${\cal A}_r$ such that a bijective
correspondence between solutions $u$ of $(1.1)$ and $\hat{u}$ of
$(1.2)$ exists, where $\hat{A}_1,...,\hat{A}_n$ are PDO over the
Cayley-Dickson algebra, $Q$ is a polynomial (or power series) with
${\cal A}_r$ coefficients, $\hat{g}$ and $\hat{u}$ are ${\cal
A}_r$-valued functions.}
\par {\bf Proof.} I. Let PDE $(1.1)$ be in a variable $x=(x_1,...,x_n)$ belonging to a domain
$U$ in the Euclidean space ${\bf R}^n$, where each variable
$x_1$,...,$x_n$ belongs to the real field $\bf R$, where $U$ is open
in $V$ by the conditions of this theorem. The subset $V$ in ${\bf
R}^n$ is canonically closed, which means by the definition that the
closure $cl (Int (V))$ of the interior $Int (V)$ of $V$ coincides
with $V$. To this variable $x=(x_1,...,x_n)\in U$ we pose a variable
$z=z(x)$ by the formula
\par $(1.3)$ $z=x_1i_{l_1}+...+x_ni_{l_n}\in V$ with $l_s\ne l_p$ for each
$s\ne p$, where $l_1,...,l_n$ are fixed nonnegative integers, $V$
notates the corresponding domain in the Cayley-Dickson algebra
${\cal A}_t$ with $n\le 2^t$, $2\le t$. The family $ \{ i_0,
i_1,..., i_{2^r-1} \} $ denotes the standard basis of the
Cayley-Dickson algebra ${\cal A}_r$ over $\bf R$ such that $i_0=1$,
$i_j^2=-1$, $i_ji_k=-i_ki_j$ for each $j\ge 1$ and $k\ge 1$ with
$j\ne k$. Thus to each $x\in U$ a unique $z=z(x)$ is posed so that
$V=\{ z\in {\cal A}_t: ~ z=z(x), x\in U \} $. Vise versa to each
$z\in V$ a unique $x\in U$ corresponds, \par $(1.4)$ $x_j=\pi
_{l_j}(z)$ for each $j$, where $\pi _j : {\cal A}_t\to \bf R$ is an
$\bf R$-linear operator prescribed by the formulas:
\par $(1.5)$ $\pi _j(z)=z_j=(-zi_j+ i_j(2^t-2)^{-1} \{ -z
+\sum_{k=1}^{2^t-1}i_k(zi_k^*) \} )/2$ for each $j=1,2,...,2^t-1$,
\par $(1.6)$ $\pi _0(z)=z_0=(z+ (2^t-2)^{-1} \{ -z +
\sum_{k=1}^{2^t-1}i_k(zi_k^*) \} )/2$, \\
where $2\le t\in \bf N$, $z$ is a Cayley-Dickson number in ${\cal
A}_t$ presented as
\par $(1.7)$ $z=z_0i_0+z_1i_1+...+z_{2^t-1}i_{2^t-1}\in {\cal
A}_t$, $z_j\in \bf R$ for each $j$, $i_k^* = {\tilde i}_k = - i_k$
for each $k>0$, $i_0=1$, $z^*=z_0i_0-z_1i_1-...-z_{2^t-1}i_{2^t-1}$
(see Formulas II$(1.1)-(1.3)$ in \cite{ludifeqcdla}).
\par Thus to each basic vector $e_j=(0,...,0,1,0,...)$ in ${\bf
R}^n$ with $1$ at $j$-th place the basic generator $i_{l_j}$ in the
Cayley-Dickson algebra ${\cal A}_t$ is counterposed according to the
mapping ${\bf R}^n\ni x\mapsto z(x)\in {\cal A}_t$. Therefore
$(x,y)=Re (z(x)z^*(y))$ for each $x$ and $y$ in ${\bf R}^n$, where
$(x,y)=\sum_{j=1}^nx_jy_j$ is the scalar product in the Euclidean
space ${\bf R}^n$, $~ Re(w)=(w+w^*)/2$ is the real part of $w$ for
each $w\in {\cal A}_t$. Particularly, for $n=3$ the vector product
$x\times y$ can be expressed as $Im (z(x)z(y))$ with $1\le l_j$ for
each $j=1,2,3$ and with $i_{l_1}i_{l_2}=i_{l_3}$, for example,
$l_1=1$, $l_2=2$, $l_3=3$, where $Im (w) =w-Re (w)$ denotes the
imaginary part of a Cayley-Dickson number $w\in {\cal A}_t$.
\par  II. To each function $f: U\to {\bf  R}$ a unique function
$h^f(z)$ corresponds such that $h: V\to {\bf R}$ and
\par $(1.8)$ $h^f(z(x))=f(x)$ for each $x\in U$.
\par For sufficiently times differentiable function $f$ to each
partial derivative $\partial f(x)/\partial x_j$ we pose $\partial
h^f(z)/\partial z_{l_j}$ due to Formulas $(1.4)-(1.7)$ and so on by
induction to each PDO $Af(x)=\sum_{\alpha } c_{\alpha }(x)\partial
^{|\alpha |}f(x)/\partial x_1^{\alpha _1}...x_n^{\alpha _n}$ the PDO
$\hat{A}h^f(z)=\sum_{\alpha } h^{c_{\alpha }}(z)\partial ^{|\alpha
|}h^f(z)/\partial z_{l_1}^{\alpha _1}...z_{l_n}^{\alpha _n}$
corresponds with the help of Formula $(1.8)$, where $c_{\alpha }(x)$
are coefficients, $\alpha = (\alpha _1,...,\alpha _n)$, $|\alpha
|=\alpha _1+...+\alpha _n$, $\alpha _j$ is a nonnegative integer for
each $j$.

\par III. To a function $g(x)$ with values in ${\bf R}^k$ we pose a function
${\hat g}(z(x))=g_1(x)i_0+...+g_k(x)i_{k-1}$ with values in ${\cal
A}_{t_1}$, where $2^{t_1-1}<k\le 2^{t_1}$. Then to $u$ we pose a
function ${\hat u}(z(x))=u_1(x)i_{q_1}+...+u_m(x)i_{q_m}$ having
values in ${\cal A}_{t_2}$ with fixed nonnegative integers
$q_1,...,q_m$ such that $q_s\ne q_p$ for each $s\ne p$. Next we
choose $r\ge \max (t, t_1, t_2, 2)$, hence ${\cal
A}_t\hookrightarrow {\cal A}_r$ and similar embeddings are for $t_1$
and $t_2$ instead of $t$ also. For example, $q_j=j-1$ and $l_j=j-1$
can be taken for each natural number $j=1, 2,...$. \par Then
particularly \par $\sum_j\partial {\hat u}_{l_j}(z)/\partial
z_{l_j}=Re (\sigma \hat{u}^*(z))$ corresponds to \par $div ~
u(x)=(\nabla ,u)(x)= \sum_j \partial u_j(x)/\partial x_j$, where
\par $\sigma  \hat{f}(z) = \sum_j (\partial {\hat f}(z)/\partial
z_{l_j})i_{q_j}$; also \par $\sigma \hat{u}_s(z)$ to \par $grad ~
u_s(x) =\nabla ~ u_s(x)=\sum_j(\partial u_s(x)/\partial x_j)e_j$. In
particular, for $n=3$ the operator $- Im (\sigma \hat{u}(z))$ to
$rot ~ u(x)=\nabla \times u(x)$ corresponds with $q_1=1$, $q_2=2$,
$q_3=3$ or more generally with natural numbers $q_j\ge 1$ such that
$i_{q_1}i_{q_2}=i_{q_3}$.
\par Therefore taking $Q(\hat{A}_1,...,\hat{A}_v,\hat{u})={\hat
P}(\hat{A}_1,...,\hat{A}_v,(\pi _{q_1}\hat{u},...,\pi
_{q_m}\hat{u}))$ we get PDE $(1.2)$ instead of $(1.1)$, since the
Cayley-Dickson algebra ${\cal A}_r$ is power associative: $z^{\phi
}z^{\psi }=z^{\phi +\psi }$ for each natural numbers $\phi $ and
$\psi $, where an order of multiplications in $Q$ over ${\cal A}_r$
is essential (see also Section 2 in \cite{ludfov}). The
correspondences described above between domains, functions and
classes of differentiable functions are bijective. Moreover, the
mappings $U\ni x\mapsto z(x)\in V$, $g\mapsto \hat{g}$, $P\mapsto
\hat{P}$, $f\mapsto h^f$, $u\mapsto \hat{u}$ are bijective
isometries.  Thus PDE $(1.1)$ and $(1.2)$ are equivalent: to each
solution of $(1.1)$ whenever it exists a unique solution of $(1.2)$
corresponds and vise versa, since $P(A_1,...,A_v,u)$ and
$Q(\hat{A}_1,...,\hat{A}_v,\hat{u})$ exist and converge
simultaneously.

\par {\bf 2. Theorem.} {\it Suppose that conditions of Theorem 1 are
fulfilled with $P$ being a polynomial in PDOs $A_1,...,A_v$, suppose
also that each PDO $A_s$ is a linear combination of elliptic
operators of finite orders. Then Dirac-type operators $\Upsilon _j$,
$j=1,2,...,l$ exist such that the principal symbol of the PDO
$Q(\hat{A}_1,...,\hat{A}_v,\cdot )$ is $S(\Upsilon _1,...,\Upsilon
_l,\cdot )$, where $S$ is a polynomial in $\Upsilon _1,...,\Upsilon
_l$ over the Cayley-Dickson algebra ${\cal A}_p$ with $p\ge r\ge
2$.}
\par {\bf Proof.} In view of Theorem 1 PDE $(1.1)$ is equivalent to
$(1.2)$. Applying Theorem 2.1 in \cite{ludrimut2014} we get a
decomposition of each PDO $\hat{A}_s$ with the help of Dirac-type
operators $\Upsilon _j$, $j=1,2,...,l$, where $l$ is a natural
number. Each PDO $\hat{A}_s$ has a decomposition over the
Cayley-Dickson algebra ${\cal A}_{t_s}$ with $t_s\ge r\ge 2$. Taking
$p=\max (t_s: s=1,...,v)$ we get the decomposition over ${\cal A}_p$
for all PDOs $\hat{A}_1,...,\hat{A}_v$. Therefore, the principal
symbol of the PDO $Q(\hat{A}_1,...,\hat{A}_v,\hat{u})$ acting on a
function $\hat{u}$ becomes a polynomial $S(\Upsilon _1,...,\Upsilon
_l,{\hat u})$ in $\Upsilon _j$ with $j=1,2,...,l$.

\par {\bf 3. Proposition.} {\it If conditions of Theorem 1 are
satisfied and each operator $A_s$ is of the second order with
constant coefficients, $s=1,...,v$, then Dirac-type operators
$\Upsilon _j$ with $j=1,2,...,l$ exist such that
\par $Q(\hat{A}_1,...,\hat{A}_v,\hat{u})=S(\Upsilon _1,...,\Upsilon
_l,\hat{u})$, \\ where $S$ is a polynomial (or power series) in
$\Upsilon _1,...,\Upsilon _l$ and $\hat{u}$ over the Cayley-Dickson
algebra ${\cal A}_p$ with $p\ge r\ge 2$.}
\par {\bf Proof.} By virtue of Theorem 2.1 and Example 2.6 in \cite{ludrimut2014}
each PDO $\hat{A}_s$ of the second order with constant coefficients
is a polynomial of the second order in a finite number of Dirac-type
operators $\Upsilon _j$ over the Cayley-Dickson algebra ${\cal
A}_{t_s}$ with $t_s\ge r\ge 2$. Taking $Q$ given by Theorem 1 and
substituting PDOs $\hat{A}_1,...,\hat{A}_v$ with polynomials of the
second order in Dirac-type operators we deduce that
$Q(\hat{A}_1,...,\hat{A}_v,\hat{u})=S(\Upsilon _1,...,\Upsilon
_l,\hat{u})$, where $S$ is a polynomial (or power series) in
$\Upsilon _1,...,\Upsilon _l$ and $\hat{u}$ over the Cayley-Dickson
algebra ${\cal A}_p$ with $p=\max (t_s: s=1,...,l)$.

\section{Nonlinear Sobolev-Burgers PDE.}
\par {\bf 4. Sobolev-Burgers PDEs.}
A more general Sobolev type PDE is with $Q(\frac{\partial }{\partial
t})$ instead of $\frac{\partial }{\partial t}$ in the Burgers PDE,
where $Q(t)=t^m+c_{m-1}t^{m-1}+...+c_0$ is a polynomial with real or
complex coefficients $c_0,...,c_{m-1}$, where $m\ge 1$ is a natural
number. We take them in the form:
$$(4.1)\quad Q(\frac{\partial }{\partial t}) (- \Delta _x^2 + \alpha \Delta _x + \beta I)
u(t,x) + \gamma \frac{\partial u^2(t,x)}{\partial x_1} + \varsigma
u^2(t,x)=0,$$ where $\alpha \ne 0$, $ ~ \beta $, $\gamma $ and
$\varsigma $ are real or complex constants, $|\gamma |+|\varsigma
|>0$, $ ~ t\ge 0$, $x\in {\bf R}^n$, $x=(x_1,..,x_n)$, $$\Delta _x
=\sum_{j=1}^n \frac{\partial ^2}{\partial x_j^2}$$ denotes the
Laplace operator, $n\ge 2$, $~I$ is the unit operator (see Ch. 3,
Sect. 6 in \cite{svalkorplb}). Classically real constants $\alpha >
0$, $\beta $, $\gamma
>0$ and $\varsigma $ are considered. In the latter case by changing the
variables $t$ and $x$ and with the help of the dilation $u\mapsto
bu$, where $b$ is a real constant, one can choose $\alpha =1$ and
$\gamma =1/2$. Particularly if $\beta =0$ and $\varsigma =0$ and
$m=1$, this corresponds to the Burgers PDE.
\par It is useful to take the complexified Cayley-Dickson algebra
${\cal A}_{r,C} = {\cal A}_{r} \oplus ({\cal A}_{r}{\bf i})$, where
${\bf i}^2=-1$, $~ {\bf i}b=b{\bf i}$ for each $b\in {\cal A}_r$,
$~2\le r<\infty $. That is any complexified Cayley-Dickson number
$z\in {\cal A}_{r,C}$ has the form $z=x+{\bf i}y$ with $x$ and $y$
in ${\cal A}_r$, $x=x_0i_0+x_1i_1+...+x_{2^r-1}i_{2^r-1}$, while
$x_0,...,x_{2^r-1}$ are in $\bf R$ (see also the notation in
Subsection 1). The real part of $z$ is $Re (z)=x_0=(z+z^*)/2$, the
imaginary part of $z$ is defined as $Im (z)=z-Re (z)$, where the
conjugate of $z$ is $z^*=\tilde{z} = Re (z) - Im (z)$, that is
$z^*=x^*-{\bf i}y$ with $x^*=x_0i_0-x_1i_1-...-x_{2^r-1}i_{2^r-1}$.
Then put $|z|^2=|x|^2+|y|^2$, where
$|x|^2=xx^*=x_0^2+...+x_{2^r-1}^2$. Mention that the operator $\pi
_j$ (see Formulas $(1.5), (1.6)$ above) has the natural $\bf
C$-linear extension: $\pi _j(av+bz)=a\pi _j(v)+b\pi _j(z)$ for every
$a$ and $b$ in ${\bf C}$ and $v$ and $z$ in ${\cal A}_{r,C}$. In
another words the operator $\pi _j$ is $\bf C$-homogeneous and
${\cal A}_{r,C}$-additive.
\par Using Theorem 1 of the preceding section we write the
corresponding PDE over the Cayley-Dickson algebra ${\cal A}_r$ with
$2\le r$ and $2^r>n$. Put $x_j=z_j$ for each $j\ge 1$, where
$z_0i_0+z_1i_1+z_2i_2+... =z\in {\cal A}_r$, then $$\Delta _x= -
\sigma _z^2, \mbox{ where  }\sigma
_zf(z)=\sum_{j=1}^ni_j^*\frac{\partial f(z)}{\partial z_j},$$
$\frac{\partial }{\partial x_1} = \pi _1(\sigma _z)$, where $\pi _j$
is provided by Formula $(1.5)$ for $j\ge 1$. Functions $u$ and $f$
are supposed to be sufficient times differentiable by the
corresponding variables, so that $u$ and $f$ more generally may have
values not only in $\bf R$, but in ${\bf C}$ also.
\par  Let $\Omega $ be a set supplied
with an algebra $\cal F$ of its subsets and let $\mu : {\cal F}\to
[0,1]$ be a probability measure and let $\cal F$ be $\mu $-complete.
We take the PDE:
$$(4.2)\quad E\{ Q(\frac{\partial }{\partial t})
(- \Delta _x^2 + \alpha \Delta _x + \beta I)u(t,x;\omega ) +  \gamma
\frac{\partial u^2(t,x;\omega )}{\partial x_1} + \varsigma
u^2(t,x;\omega ) \} =0,$$ where $\omega \in \Omega $, $u(t,x;\omega
)$ is a random function, $Eg$ denotes a mean value (expectation) of
a random variable $g$ whenever it exists:
$$(4.3)\quad Eg = \int_{\Omega } g(\omega )\mu (d\omega ).$$
In order to use the noncommutative integral over ${\cal A}_r$
approach the following PDE generalizing $(4.2)$ is written in the
form:
$$(4.4)\quad E \{ Q(\frac{\partial }{\partial t}) {\bf S}_0 u(t,x,y; \omega )
+ \gamma \pi _1(\sigma _x+\sigma _y)(u^2(t,x,y;\omega ) ) +
\varsigma u^2(t,x,y;\omega ) \} |_{x=y} =0,$$ where $x$ and $y$ are
in $V\subset {\cal A}_r$, $V={\bf R}i_1\oplus ...\oplus {\bf R}i_n$,
$$(4.5)\quad {\bf S}_0 = - (\sigma _x^2+ \sigma _y^2)^2+ a(\sigma _x^2+\sigma _y^2) + bI,$$
$a\in {\bf C}\setminus \{ 0 \} $, $~b\in {\bf C}$, the Dirac
operator in $(4.5)$ is more general:
$$(4.6)\quad \sigma _xf(x)= \sum_{j=0}^{2^r-1} i_j^* (\partial
f(x)/\partial x_{\xi (j)}) {\psi }_j$$ with real constants $\psi
_j\in \bf R$ so that $\psi _0^2+...+\psi_{2^r-1}^2>0$. Particularly
it is convenient to choose $\psi _j=2^{-1/2}$ for each $j=1,...,n$;
$\psi _j=0$ otherwise for $j=0$ or $j>n$; $\xi (j)=j$ for each
$j=0,...,2^r-1$. In this particular case $-\frac{1}{2}\Delta
_x=\sigma _x^2$.

\par {\bf 5. Auxiliary PDE.} Consider the auxiliary PDE
$$(5.1)\quad \{ {\bf S}_{2,a} v(x,y)  +
q_1 \pi _1(\sigma _x+\sigma _y)(v^2(x,y)) + q_2v^2(x,y) \} |_{x=y}
=0,\mbox{ where}$$
\par $(5.2)$ ${\bf S}_{2,a} = a_1 (\sigma _x^2+ \sigma _y^2)^2+ a_2(\sigma _x^2+\sigma _y^2) + a_3I,$
\\ $a_1\in {\bf C}\setminus \{ 0  \} $, $ ~a_2\in {\bf C}$ and $a_3\in {\bf C}$ are constants, $a_1$ is nonzero,
$a=(a_1,a_2,a_3)$, $ ~ {\bf C}={\bf R}\oplus {\bf R}{\bf i}$; $ ~
q_1=-2a_1p_1$, $ ~ q_2= - 2 a_1p_2$,  $~p_1$ and $p_2$ are in ${\bf
C}$ with $|p_1|+|p_2|>0$; where $i_j{\bf i}={\bf i}i_j$ for each
$j$, ${\bf i}^2=-1$. Remind that the family $ \{ i_0, i_1,...,
i_{2^r-1} \} $ denotes the standard basis of the Cayley-Dickson
algebra ${\cal A}_r$ over $\bf R$ such that $i_0=1$, $i_j^2=-1$,
$i_ji_k=-i_ki_j$ for each $j\ge 1$ and $k\ge 1$ with $j\ne k$. For
its solution take two PDOs:
\par $(5.3)$ ${\bf S}_1=\sigma _x^2-\sigma_y^2$ \par and
${\bf S}_2= {\bf S}_{2,a}$, where $a=(a_1,a_2,a_3)$. Suppose that a
function $F(x,y)=F_{\bf a}(x,y)$ satisfies the conditions
\par $(5.4)$ ${\bf S}_jF(x,y)=0$ for $j=1$ and $j=2$. Put also \par $(5.5) \quad K(x,y)=F(x,y)+{\sf
A}K(x,y)$ with
$$(5.6)\quad {\sf A}K(x,y)= p_1\pi _1(\mbox{
}_{\sigma }\int_{w_0} ^x (\mbox{ }_{\sigma }\int_{w_0}^y (\mbox{
}_{\sigma }\int_w^{\infty } F(z,v)K(w,z)dz)dv)dw)$$
$$+p_2 \mbox{
}_{\sigma }\int_{w_0} ^y(\mbox{ }_{\sigma }\int_{w_0} ^x (\mbox{
}_{\sigma }\int_{w_0}^s (\mbox{ }_{\sigma }\int_w^{\infty }
F(z,v)K(w,z)dz)dv)dw)ds,$$ where $p_1$ and $p_2$ are real or complex
parameters with $|p_1|+|p_2|>0$; $ ~ K=K_{{\bf a},{\bf q}}$, $ ~
{\bf a}=(a_1,a_2,a_3),$ $ ~ {\bf q}=(q_1,q_2)$, $~ w_0$ is marked
point in a domain ${\cal U}$ satisfying Conditions 3.1$(D1)$ and
$(D2)$ in \cite{lucveleq2013}.
\par According to Formulas $(5.3)-(5.6)$ it is convenient to consider the functions $F$ and $K$ of the form:
\par $(5.7)$ $F(x,y)=F(\frac{x+y}{2})$ and $K(x,y)=K(\frac{x+y}{2})$.
\par $(5.8)$ Let $F$ and $K$ be with values in $\bf R$ or $\bf C$ and
let they satisfy other Conditions of Proposition 3.4 in
\cite{lucveleq2013} with $m=4$. \par Apparently any pair of
operators from the family $ \{ {\bf S}_1, {\bf S}_{2,a}, bI, \pi _j:
b\in {\bf C}; j=0,..., 2^{r-1}\}$ commutes, since ${\bf S}_1$ and
${\bf S}_{2,a}$ are PDOs with constant real or complex coefficients,
where $I$ notates the unit operator. From this proposition and
Formula $(3.2)$ in \cite{lucveleq2013} we infer that
$$ ({\bf S}_{2,a}-a_3I)K(x,y)=({\bf S}_{2,a}-a_3I)F(x,y)+ $$ $$p_1 \pi _1 ({\bf S}_{2,a}-a_3I)((\mbox{ }_{\sigma }\int_{w_0} ^x
(\mbox{ }_{\sigma }\int_{w_0}^y (\mbox{ }_{\sigma }\int_w^{\infty }
F(z,v)K(w,z)dz)dv)dw))$$
$$ +p_2 ({\bf S}_{2,a}-a_3I)(\mbox{
}_{\sigma }\int_{w_0} ^y(\mbox{ }_{\sigma }\int_{w_0} ^x (\mbox{
}_{\sigma }\int_{w_0}^s (\mbox{ }_{\sigma }\int_w^{\infty }
F(z,v)K(w,z)dz)dv)dw)ds) $$
$$= -a_3F(x,y)+ p_1\pi _1 \{ a_1 (\sigma _x^2+ \sigma _y^2)+ a_2I \} (\mbox{ }_{\sigma }\int_x^{\infty }
F(z,y)K(x,z)dz)$$ $$+ p_2 a_1 (\sigma _x+ \sigma _y) (\mbox{
}_{\sigma }\int_x^{\infty } F(z,y)K(x,z)dz)$$
$$+ p_2 a_2(\mbox{ }_{\sigma
}\int_{w_0} ^y (\mbox{ }_{\sigma }\int_x^{\infty }
F(z,v)K(x,z)dz)dv) $$
$$ = -a_3F(x,y)+ p_1\pi _1 \{ a_1 (\mbox{}^2\sigma_x^2+\mbox{ }^2\sigma
_z^2)+a_2I \}  (\mbox{ }_{\sigma } \int_x^{\infty } F(z,y)K(x,z)dz)
$$
$$+p_1\pi _1 a_1(A_2(F,K)(x,y)+B_2(F,K)(x,y))$$
$$+p_2 a_1 (\mbox{}^2\sigma_x + \mbox{ }^2\sigma
_z)  (\mbox{ }_{\sigma } \int_x^{\infty } F(z,y)K(x,z)dz)
+p_2a_1(A_1(F,K)(x,y)+B_1(F,K)(x,y))$$
$$+p_2 a_2(\mbox{ }_{\sigma
}\int_{w_0} ^y (\mbox{ }_{\sigma }\int_x^{\infty }
F(z,v)K(x,z)dz)dv)  .$$
 Therefore
$$(5.9)\quad {\bf S}_{2,a}K(x,y)= \mbox{ }^2{\bf S}_{2,a}p_1\pi _1((\mbox{ }_{\sigma }\int_{w_0} ^x (\mbox{ }_{\sigma }\int_{w_0}^y
(\mbox{ }_{\sigma }\int_w^{\infty } F(z,v)K(w,z)dz)dv)dw))$$ $$ +
p_1\pi _1a_1(A_2(F,K)(x,y)+B_2(F,K)(x,y))  +
p_2a_1(A_1(F,K)(x,y)+B_1(F,K)(x,y)) .$$  By virtue of Corollary 3.6
in \cite{lucveleq2013}
$$(5.10)\quad A_2(F,K)(x,y)+B_2(F,K)(x,y)=$$
$$ - \sigma _x [F(x,y)K(x,x)] - ~\mbox{}^2\sigma
_x[F(z,y)K(x,z)]|_{z=x}$$ $$ - \mbox{}^1\sigma
_z[F(z,y)K(x,z)]|_{z=x} + ~\mbox{}^2\sigma _z[F(z,y)K(x,z)]|_{z=x}
\mbox{ and}$$  \par $(5.11)$ $A_1(F,K)(x,y)+B_1(F,K)(x,y)= - 2
F(x,y)K(x,x)$. \par Moreover, this operator ${\sf A}$ restricted on
any compact domain $V$ in ${\cal U}$ is compact due to Proposition
4.1 in \cite{lucveleq2013} and Formula $(5.6)$. The PDOs ${\bf S}_j$
commute with $\sigma _x$ and $\sigma _y$ for $j=1$ and $j=2$, since
${\bf i}i_k=i_k{\bf i}$ for each $k$. Thus from Formulas $(3.2)$ in
\cite{lucveleq2013} and $(5.9)$ and $(5.10)$ given above the theorem
follows.

\par {\bf 6. Theorem.} {\it Let suppositions $(5.8)$ be
fulfilled, then the PDE $(5.1)$ has a solution $K_{{\bf a},{\bf q}}$
provided by Formulas $(5.4)-(5.6)$ on $({\cal U}\cap V)^2$ for a
sufficiently small $0<|p_1|+|p_2|$.}

\par Below a random operator valued measure approach
to a solution of Sobolev-Burgers PDEs is described. To avoid a
misunderstanding we first recall necessary definitions and describe
a notation and provide necessary statements.

\par {\bf 7. Definition. Orthogonal random operator valued measure.}
\par For Banach spaces $X$ and $Y$ both over $\bf F$ by
${\cal L}(X,Y)$ is denoted the space of all bounded linear operators
$J$ from $X$ into $Y$, where either ${\bf F}=\bf C$ or ${\bf F}=\bf
R$. If it is supplied with the operator norm topology $\tau _{|\cdot
|}$, then $({\cal L}(X,Y), \tau _{|\cdot |})$ is a Banach space,
where $|J|=\sup_{x\in X, ~ 0<|x|\le 1 ~} |Jx|_Y/|x|_X$ for any $J\in
{\cal L}(X,Y)$, $~|\cdot |_X$ denotes a norm on $X$. A strong
operator topology $\tau _s$ on ${\cal L}(X,Y)$ possesses a base of
neighborhoods of $0$ of the form $Z_s(\epsilon , x) := \{ J\in {\cal
L}(X,Y): ~ |Jx|_Y<\epsilon \} $, where $x\in X$, $~0<\epsilon $. A
weak operator topology $\tau _w$ is induced by a base of
neighborhoods of zero $Z_w(\epsilon , x,y') := \{ J\in {\cal
L}(X,Y): ~ |y'(Jx)|<\epsilon \} $, where $x\in X$, $ ~ y'\in Y'$, $
~ 0<\epsilon $, $ ~ Y'$ is a topological dual space of all
continuous linear functionals $y': Y\to {\bf F}$. For short the
topological locally convex vector space $({\cal L}(X,Y), \tau
_{\kappa })$ will also be denoted by ${\cal L}(X,Y)_{\kappa }$,
where $\kappa \in \{ {|\cdot |}, s , w \} $.
\par Let $(\Omega ,{\cal F}, \nu)$ be a measure space, where $\nu$ is a
nonnegative $\sigma $-finite and $\sigma $-additive measure on a
$\sigma $-algebra $\cal F$ of a set $\Omega $. For the Banach space
$Y$ over $\bf F$ we notate a Banach space of all $({\cal F},
B(Y))$-measurable functions $f: \Omega \to Y$ such that
$$(7.1)\quad | f |_v := [\int_{\Omega } | f(\omega ) |^v_Y \nu
(d\omega )]^{1/v}<\infty $$ by $L^v(\Omega, {\cal F}, \nu, Y)$,
where $1\le v<\infty $, while $B(Y)$ is the Borel $\sigma $-algebra
on $Y$.
\par Then $L^v(\Omega ,{\cal F}, \nu; {\cal L}(X,Y)_{\kappa })$
denotes a linear topological space of all $({\cal F}, B({\cal
L}(X,Y)_{\kappa }))$ measurable functions $G: \Omega \to {\cal
L}(X,Y)_{\kappa }$ such that either for $\kappa = {|\cdot |}$
$$(7.2)\quad |G|_{v} := [\int_{\Omega } |
G(\omega )|^v \nu(d\omega )]^{1/v}<\infty $$ or for $\kappa =s$ and
each $x\in X$
$$(7.3)\quad \rho (G)_{v,x} := [\int_{\Omega } |
G(\omega )x|^v_Y \nu(d\omega )]^{1/v}<\infty $$ or for $\kappa =w$
and every $x\in X$ and $y'\in Y'$
$$(7.4)\quad \rho (G)_{v,x,y'} := [\int_{\Omega } |y'(
G(\omega )x)|^v \nu(d\omega )]^{1/v}<\infty .$$  \par Let  $\Lambda
$ be a set and $\cal M$ be a semi-ring of subsets in it. Suppose
that $X$ and $Y$ are Hilbert spaces over $\bf F$, also $(\Omega
,{\cal F}, \mu )$ is a probability space (see Subsection 4) and to
each $M\in \cal M$ a ${\cal L}(X,Y)_{\kappa }$ valued random
operator $H(M)(\omega )$ corresponds and fulfilling Conditions
$(7.5)-(7.7)$:
\par $(7.5)$ for each $M\in \cal M$  \par
$\mu \{ \omega \in \Omega : {\cal D}(H(M)(\omega ))=X \} =1$ and
$H(M)\in L^2(\Omega ,{\cal F}, \mu ; {\cal L}(X,Y)_{\kappa })$ and
$\mu \{ \omega \in \Omega : H(\emptyset )(\omega ) =0 \} =1$, where
\\  ${\cal D}(H(M)(\omega ))$
is a linear definition domain of $H(M)(\omega )$ dense in the
Hilbert space $X$; $ ~ H(M)(\omega ): {\cal D}(H(M)(\omega ))\to Y$;
shortly $H(M)$ is written instead of $H(M)(\omega )$;
\par $(7.6)$ for every $M_1$ and $M_2$ in $\cal M$ with $M_1\cap M_2=\emptyset $
\par $H(M_1\cup M_2) = H(M_1)+H(M_2) ~ (mod ~ \mu )$;
\par $(7.7)$ for every $M_1$ and $M_2$ in $\cal M$ the mean value
\par $\hat{{\sf m}}(M_1,M_2):= E (H^*(M_1)(\omega )H(M_2)(\omega ))$ exists
and belongs to ${\cal L}(X,X)_w$, where $\kappa \in \{ {|\cdot |}, s
\}$; $ ~ H^*(M)$ is an adjoint operator, that is
$(H^*(M)y,x)_X=(y,H(M)x)_Y$ for each $x\in {\cal D}(H(M))$ and $y\in
{\cal D}(H^*(M))$, where ${\cal D}(H^*(M))\subset Y$, since $Y$ is
the Hilbert space and hence $Y'$ is isomorphic with $Y$, $ ~ H^*(M):
{\cal D}(H^*(M))\to X$. Shortly $\hat{{\sf m}}(M)$ also will be
written instead of $\hat{{\sf m}}(M,M)$.
\par $(7.8)$ $\hat{{\sf m}}(M_1,M_2)=0$ if $M_1\cap M_2=\emptyset $.
\par The family of random operators $\{ H(M): ~ M\in {\cal M} \} $
satisfying Conditions $(7.5)-(7.7)$ is called an elementary random
operator valued measure and $\hat{{\sf m}}$ is called its structural
function. It is called orthogonal if $(7.8)$ also is valid.
\par For a Hilbert space $X$ over $\bf C$ by
$(x,y)_X=(x,y)$ is denoted a scalar product on $X$ with values in
$\bf C$ so that $(\beta x,y)=\beta (x,y)$ and $(x,\beta
y)=\bar{\beta }(x,y)$ and $(x,y+u)=(x,y)+(x,u)$ and
$(x,y)=\overline{(y,x)}$ for every $x$ and $y$ and $u$ in $X$ and
$\beta \in {\bf C}$, where $\bar{\beta }$ is the complex conjugated
number of $\beta $. An induced norm is $|x|_X=\sqrt{(x,x)_X}$ for
each $x\in X$. Particularly, if $X$ is a Hilbert space over $\bf R$,
then $(x,\beta y)=\beta (x,y)$ and $(x,y)=(y,x)\in \bf R$, where
$\beta \in \bf R$.

\par {\bf 8. Lemma.} {\it Suppose that $H$ is an elementary orthogonal random
operator valued measure with a structural function $\hat{{\sf m}}$.
If $M_1$ and $M_2$ are in $\cal M$, then
\par $(8.1)$ $\hat{{\sf m}}(M_1,M_2)=\hat{{\sf m}}(M_1\cap M_2)$.
\par For each $M\in \cal M$ an operator $\hat{{\sf m}}(M)$ is nonnegative
definite: \par $(8.2)$ for every $c_1,...,c_n\in \bf C$ and
$x_1,...,x_n\in X$  $$ \sum_{j,k=1}^n c_j\bar{c}_k (\hat{\sf
m}(M)x_j,x_k)_X\ge 0.$$ If $M_1\cap M_2=\emptyset $, then \par
$(8.3)$ $\hat{{\sf m}}(M_1\cup M_2)=\hat{{\sf m}}(M_1)+ \hat{{\sf
m}}(M_2)$.
\par For every $x$ and $y$ in $X$ and $N\in {\cal M}$
\par $(8.4)$ $|\hat{{\sf m}}_{x,y}(N)|\le \sqrt{\hat{{\sf m}}_{x,x}(N)\cdot \hat{\sf
m}_{y,y}(N)}$, \\ where $\hat{{\sf m}}_{x,y}(N) := (\hat{\sf
m}(N)x,y)_X$.}

\par {\bf Proof.} Assertions $(8.1)$ and $(8.2)$ follow from $(7.5)-(7.8)$, since
$M_1\setminus M_2\in \cal M$ and $M_2\setminus M_1\in \cal M$,
$M_1\cap (M_2\setminus M_1)=\emptyset $ and $M_2\cap (M_1\setminus
M_2)=\emptyset $ and $\hat{{\sf m}}(M_1,M_2)=EH^*(M_1\cap
M_2)H(M_1\cap M_2) + EH^*(M_1\setminus M_2) H(M_1\cap M_2) +
EH^*(M_1\cap M_2)H(M_2\setminus M_1)+ EH^*(M_1\setminus
M_2)H(M_2\setminus M_1)=\hat{{\sf m}}(M_1\cap M_2)$.\par Therefore,
if $M_1\cap M_2=\emptyset $, then
\par $\hat{{\sf m}}(M_1\cup M_2)=EH^*(M_1)H(M_1) + EH^*(M_1) H(M_2) +
EH^*(M_2)H(M_1)+ EH^*(M_2)H(M_2)=\hat{{\sf m}}(M_1)+\hat{{\sf
m}}(M_2)$. Using $(7.5)$ and $(7.7)$ one gets that
$$ \sum_{j,k=1}^n c_j\bar{c}_k (\hat{{\sf m}}(M)x_j,x_k)= E(H(M)x,H(M)x)\ge 0,$$
where $x=c_1x_1+....+c_nx_n$. Then from \par $|\hat{\sf
m}_{x,y}(N)|=|E(H(N)x,H(N)y))|\le E(|H(N)x|\cdot |H(N)y|)$\par $\le
\sqrt{(E|H(N)x|^2)\cdot (E|H(N)x|^2)}$, \\ Inequality $(8.4)$
follows.

\par {\bf 9. Definition.}
Let $L^0({\cal M},X)$ denote a linear space of all step functions
$f: \Lambda \to X$ such that $f(\lambda ) =\sum_{k=1}^{\iota } {\bf
\chi }_{M_k}(\lambda )a_k $, where $X$ is a Hilbert space over $\bf
F$, $ ~ M_k\in \cal M$, $ ~ a_k\in X$ for each $k=1,...,{\iota }$; $
~ {\iota }\in \bf N$; $ ~{\bf \chi }_{M}(\lambda ) $ is the
characteristic function of a set $M$, that is ${\bf \chi
}_{M}(\lambda ) =1$ for each $\lambda \in M$, whilst ${\bf \chi
}_{M}(\lambda ) =0$ for each $\lambda \notin M$.
\par For each $f\in L^0({\cal M},X)$ the
integral relative to $H$ is defined: $$(9.1)\quad \psi (f)
=\int_{\Lambda } H(d\lambda ) f(\lambda ):= \sum_{k=1}^{\iota }
H(M_k)a_k.$$
\par By $L^0(H,Y)$ we denote the family of all
random vectors $\eta =\psi (f)$ of the form $(9.1)$.
\par A sequence of random vectors $h_{\iota }$ in $L^2(\Omega , {\cal F},
\mu , Y)$ mean square converges to $h$ if $lim_{{\iota }\to \infty}
|h-h_{\iota }|_{2}=0$ and this is denoted by $l.i.m._{{\iota }\to
\infty } h_{\iota }=h$.

\par {\bf 10. Remark.} If $f$ and $g$ are in $L^0({\cal M},X)$, then ${\iota }\in \bf N$ and
$N_k\in \cal M$ for each $k=1,...,{\iota }$ can be chosen such that
$f(\lambda ) =\sum_{k=1}^{\iota } {\bf \chi }_{N_k}(\lambda )a_k $
and $g(\lambda ) =\sum_{k=1}^{\iota } {\bf \chi }_{N_k}(\lambda )b_k
$ for each $\lambda \in \Lambda $, where $a_k$ and $b_k$ are in $X$
for each $k$. Therefore from Lemma 8 and Condition $(7.7)$ it
follows that
$$(10.1)\quad E(\int_{\Lambda } H(d\lambda ) f(\lambda ),
\int_{\Lambda } H(d\lambda ) g(\lambda ))_Y = \sum_{k=1}^{\iota }
(\hat{\sf m}(N_k)a_k,b_k)_X.$$ The space $(L(X,Y), \tau _{|\cdot
|})$ is Banach. By virtue of the Banach-Steinhaus theorem (11.6.1)
in \cite{narib} or 7.1.3 in \cite{edwardsb} $(L(X,Y), \tau _{\kappa
})$ is complete as the topological locally convex vector space also
for $\kappa =s$ and $\kappa =w$.
\par Remind that a measure $\hat{{\sf m}}$ on $\cal M$ with values in $(L(X,X), \tau _{\kappa
})$ is called $\sigma $-finite, if $\Lambda = \bigcup_{k=1}^{\infty
}M_k$, where $M_k\in \cal M$ for each $k\in \bf N$. If a finitely
additive $\sigma $-finite nonnegative definite $L(X,X)_w$ valued
measure $\hat{{\sf m}}$ satisfies the semi-additivity condition:
$$(10.2)\quad \hat{{\sf m}}(N)\le \sum_{k=1}^{\infty } \hat{{\sf m}}(N_k)$$ for every $N$ and $N_k$ in
$\cal M$ fulfilling the inclusion $N\subseteq \bigcup_{k=1}^{\infty
}N_k$, that is $$ \hat{{\sf m}}_x(N):=(\hat{{\sf m}}(N)x,x)_X\le
\sum_{k=1}^{\infty } (\hat{{\sf m}}(N_k)x,x)_X$$ for each $x\in X$,
then $\hat{{\sf m}}$ has a unique extension on a minimal $\sigma
$-algebra $\sigma {\cal M}$ generated by $\cal M$, since $\hat{{\sf
m}}_x$ is with values in $[0,\infty )$ and has a unique extension on
$\sigma {\cal M}$ for each $x\in X$ (see Theorem II.2.3 in
\cite{gihmskorb}). In this case the structural function $\hat{{\sf
m}}$ will be called a structural operator valued measure. \par
Notice that $\forall M\in \sigma \cal M$ (($\hat{{\sf m}}(M)=0$)
$\Leftrightarrow $ ($\forall x\in X$ $ ~ \hat{{\sf m}}_x(M)=0$)).
Due to the scalar product properties $\hat{{\sf
m}}_{x,y}(M):=(\hat{\sf m}(M)x,y)_X$ can be expressed through a
linear combination of $\hat{\sf m}_{z,z}=\hat{{\sf m}}_z$ for
suitable vectors $z\in \{ x, y, x\pm y \} $ over ${\bf F}={\bf R}$
or $z\in \{ x, y, x\pm y , x\pm iy \} $ over ${\bf F}={\bf C}$.
Therefore $\hat{{\sf m}}$ can be extended to a complete $\sigma
$-additive measure on ${\cal B}={\cal B}_{\hat{{\sf m}}}(\Lambda )$,
where a $\sigma $-algebra $\cal B$ is the completion of $\sigma
{\cal M}$ by $\hat{\sf m}$-null sets: $\hat{{\sf m}}(N)=0$ if
$N\subset M$ and $M\in \sigma {\cal M}$ and $\hat{{\sf m}}(M)=0$. We
denote by $L^2(\Lambda , {\cal B}, \hat{{\sf m}}, X)$ the completion
of $L^0({\cal B},X)$ relative to a norm $|f|_{2,\hat{{\sf m}}} =
\sqrt{(f,f)_{\hat{{\sf m}}}}$ induced by a scalar product
$$(10.3)\quad (f,g)_{\hat{{\sf m}}} := \int_{\Lambda } (\hat{{\sf m}}(d\lambda )f(\lambda
),g(\lambda ))_X.$$ This implies that $L^0({\cal M},X)$ is a linear
subspace in $L^2(\Lambda , {\cal B}, \hat{{\sf m}}, X)$, hence a
closure $L^2({\cal M},X)$ of $L^0({\cal M},X)$ in $L^2(\Lambda ,
{\cal B}, \hat{{\sf m}}, X)$ exists. The closure of $L^0(H,Y)$ in
$L^v(\Omega , {\cal F}, \mu , Y)$ will be denoted by $L^v(H,Y)$,
where $1\le v<\infty $.
\par Formulas $(9.1)$ and $(10.1)$ induce a linear isometry $\psi $ from
$L^2({\cal M},X)$ into $L^2(H,Y)$. Hence this integral has an
extension
$$(10.4)\quad \int_{\Lambda } H(d\lambda )f(\lambda ) := \psi (f)$$
for each $f\in L^2({\cal M},X)$.
\par From Lemma 8 and Remark 10 assertions of Theorem 11 follow.

\par {\bf 11. Theorem.} {\it If Conditions $(7.5)-(7.8)$ and
$(10.2)$ are satisfied, then \par $(11.1)$ for every $f$ and $g$ in
$L^2(\Lambda , {\cal B}, \hat{{\sf m}}, X)$ and $a$ and $b$ in $\bf
C$
$$\int_{\Lambda } H(d\lambda )(af(\lambda ) + bg(\lambda ))=
a\int_{\Lambda } H(d\lambda )f(\lambda ) + b \int_{\Lambda }
H(d\lambda )g(\lambda );$$ $(11.2)$ if a sequence $f_{\iota }$
converges to $f$ in $L^2(\Lambda , {\cal B}, \hat{{\sf m}}, X)$,
then
$$\int_{\Lambda } H(d\lambda )f(\lambda ) = l.i.m._{{\iota }\to \infty }
\int_{\Lambda } H(d\lambda )f_{\iota }(\lambda ).$$}

\par {\bf 12. Remark. Extension of an elementary orthogonal random
operator valued measure.} We put \par $(12.1)$ ${\cal B}_{0,\kappa }
:= \{ N\in {\cal B}: ~ H(N)\in L^2(\Omega ,{\cal F}, \mu , {\cal
L}(X,Y)_{\kappa })$\par $  \&  ~ \hat{{\sf m}}(N)\in L(X,X)_w \} $,
where $\kappa \in \{ {|\cdot |}, s \} $.
\par Let ${\hat H} (N)=\int_{\Lambda } H(d\lambda ){\bf \chi }_N(\lambda )$,
then
\par $(12.2)$ ${\hat H}$ is defined on ${\cal B}_{0,\kappa }$ and
\par $(12.3)$ If $N_k\in {\cal B}_{0,\kappa }$ for each $k=0, 1,2
,...$ and $N_0=\bigcup_{k=1}^{\infty }N_k$ with $N_k\cap
N_l=\emptyset $ for each $k\ne l$, then ${\hat
H}(N_0)=\sum_{k=1}^{\infty }{\hat H}(N_k)$ and this series converges
in $L^2(\Omega , {\cal F}, \mu , {\cal L}(X,Y)_s)$ according to
Theorem 11, that is for each $x\in X$ \par ${\hat
H}(N_0)x=l.i.m._{\iota \to \infty } \sum_{k=1}^{\iota }{\hat H}(N_k)x$ \\
(see Definition 9) and \par $(12.4)$ $\forall N\in {\cal
B}_{0,\kappa }$ $ ~\forall M\in {\cal B}_{0,\kappa }$ $ ~ E{\hat
H}^*(N){\hat H}(M)=\hat{{\sf m}}(N\cap M)$ and \par $(12.5)$
$\forall M\in {\cal M}$ $ ~{\hat H}(M)=H(M)$.

\par {\bf 13. Definition.} A random function ${\hat H}$ satisfying
Conditions $(12.2)-(12.4)$ is called an orthogonal random operator
valued measure.
\par From Lemma 8 and Theorem 11 and Remark 12 the theorem follows.

\par {\bf 14. Theorem.} {\it If a structural function $\hat{{\sf m}}$ of an
elementary orthogonal random operator valued measure $H$ is
semi-additive (see Formula $(10.2)$), then an orthogonal random
operator valued measure $\hat H$ extension of $H$ exists.}
\par {\bf 15. Corollary.} {\it If the conditions of Theorem 14 are
satisfied, then $L^2(H,Y)$ is isomorphic with $L^2({\hat H},Y)$ and
\par $\forall f\in L^2({\hat H},Y)$ $\exists $ $\int_{\Lambda }
H(d\lambda )f(\lambda ) = \int_{\Lambda } {\hat H}(d\lambda
)f(\lambda )$.}

\par {\bf 16. Remark.} Let $H$ be an orthogonal random operator
valued measure and let $\hat{{\sf m}}$ be its complete structural
operator valued measure. For short $H$ will be written instead of
$\hat H$. For each $N\in \cal B$ and $g\in L^2(\Lambda , {\cal B},
\hat{{\sf m}}, X)$ we put
$$(16.1)\quad \eta (N):= \int_{\Lambda } H(d\lambda ) g(\lambda
){\bf \chi }_N(\lambda ), \mbox{ consequently,}$$
$$(16.2)\quad E(\eta (N),\eta (M))_{Y'}=\int_{N\cap M} (\hat{\sf m}(d\lambda )g(\lambda ),
g(\lambda ))_X$$ for all $N$ and $M$ in $\cal B$. Then we put
$$(16.3)\quad \hat{\sf n}(N) = \int_N(\hat {\sf m}(d\lambda )g(\lambda ), g(\lambda ))_X.$$
Therefore $\eta $ is an orthogonal random vector valued measure and
$\hat{\sf n}$ is its complete structural measure such that $\eta
(N)\in L^2(\Omega ,{\cal F}, \mu ; Z)$ and $\hat{\sf n}(N)\ge 0$ for
each $N\in {\cal B}$, where either $Z=(Y', |\cdot |_{Y'})$ if $\tau
_{\kappa }=\tau _{|\cdot |}$ or $Z=(Y', \sigma (Y',Y))$ if $\tau
_{\kappa }=\tau _s$, since $Y$ is the Hilbert space over $\bf F$ and
$Z={\cal L}(Y,{\bf F})_{\kappa }$, where $\sigma (Y',Y)$ denotes the
weak topology on the topological dual space $Y'$.
\par Since $Y$ is the Hilbert space, then $Y'$ and $Y$ are
isomorphic. That is $Y'$ can be replaced on $Y$ in the notation.
\par {\bf 17. Lemma.} {\it If $f\in L^2(\Lambda , {\cal B}, \hat{\sf n},
{\bf F})$ and the conditions of Remark 16 are fulfilled, then $fg\in
L^2(\Lambda , {\cal B}, \hat{{\sf m}}, X)$ and
$$(17.1)\quad \int_{\Lambda } f(\lambda )\eta (d\lambda )=\int_{\Lambda
}H(d\lambda )f(\lambda )g(\lambda ).$$}
\par {\bf Proof.} Take a fundamental sequence of step functions $f_{\iota }$
in $L^2(\Lambda , {\cal B}, \hat{\sf n}, {\bf F})$. Then
$$E(|\int_{\Lambda } (f_k(\lambda )-f_{k+l}(\lambda
))\eta (d\lambda )|^2)= \int_{\Lambda }|f_k(\lambda
)-f_{k+l}(\lambda )|^2\hat{\sf n}(d\lambda )$$ for each pair of
natural numbers $k$ and $l$, consequently, the sequence $f_{\iota
}g_{\iota }$ is fundamental in $L^2(\Lambda , {\cal B}, \hat{{\sf
m}}, X)$, hence
$$\exists ~ \lim_{\iota \to \infty }\int_{\Lambda } f_{\iota }(\lambda ) \eta
(d\lambda ) =\lim_{\iota \to \infty }\int_{\Lambda } H(d\lambda
)f_{\iota }(\lambda ) g(\lambda ) .$$ Thus for $f=\lim_{\iota \to
\infty }f_{\iota }$ Formula $(17.1)$ follows.

\par {\bf 18. Remark.} Consider an open or a canonically closed domain
$V$ in the Euclidean space ${\bf R}^k$ or in the unitary space ${\bf
C}^k$ and let $\bf l$ be a Lebesgue measure on it. Suppose that
${\cal B}(V)$ is a $\bf l$ complete $\sigma $-algebra containing the
Borel $\sigma $-algebra $B(V)$ of $V$, $ ~ B(X)$ is the Borel
$\sigma $-algebra of $(X, |\cdot |_X)$, $ ~\sigma ({\cal B}(V)\times
{\cal B})$ denotes the minimal $\sigma $-algebra containing ${\cal
B}(V)\times {\cal B}$, $~{\cal B}$ is the complete $\sigma $-algebra
on $\Lambda $ as above. \par Let $g(\tau ,\lambda )$ be a $(\sigma
({\cal B}(V)\times {\cal B}),B(X))$ measurable function from
$V\times \Lambda $ into a Hilbert space $X$ over $\bf F$, $ ~g\in
L^2(V\times \Lambda , \sigma ({\cal B}(V)\times {\cal B}), {\bf
l}\times \hat{{\sf m}}, X)$ and such that for each marked $\tau \in
V$ the vector valued function $g(\tau , \lambda )$ considered in the
$\lambda $ variable belongs to $L^2(\Lambda , {\cal B}, \hat{{\sf
m}}, X)$, where $\lambda \in \Lambda $. If $H$ is an orthogonal
random operator valued measure and $\hat{{\sf m}}$ is its complete
structural operator valued measure, then the integral
$$(18.1)\quad \xi (\tau ) = \int_{\Lambda } H(d\lambda )g(\tau
,\lambda )$$ is defined for each $\tau \in V$ and for $\mu $-a.e.
$\omega \in \Omega $ (see Theorem 14 and Corollary 15).

\par {\bf 19. Lemma.} {\it If conditions of Remark 18 are satisfied,
then the integral in Formula $(18.1)$ as a function in $\tau $ can
be defined such that $\xi (\tau )$ will be $(\sigma ({\cal
B}(V)\times {\cal F}),B(Y))$ measurable.}
\par {\bf Proof.} If $g$ is a step function $$(19.1)\quad g(\tau , \lambda
) = \sum_{k=1}^{\iota }{\bf \chi }_{N_k}(\tau ) {\bf \chi
}_{M_k}(\lambda ) x_k$$ with $N_k\in {\cal B}(V)$, $M_k\in {\cal
B}$, $x_k\in X$ for $k=1,...,\iota $, then $$\xi (\tau
)=\sum_{k=1}^{\iota } H(M_k) {\bf \chi }_{N_k}(\tau ) x_k,$$
consequently, $\xi (\tau )=\xi (\tau )(\omega )$ is $(\sigma ({\cal
B}(V)\times {\cal F}),B(Y))$ measurable, since $Y'$ and $Y$ are
isomorphic. Choose a sequence of step functions $g_{\iota }$ of the
form $(19.1)$ so that $g_{\iota }$ converges to $g$ in $L^2(V\times
\Lambda , \sigma ({\cal B}(V)\times {\cal B}), {\bf l}\times
\hat{{\sf m}}, X)$. Let $\xi _{\iota }$ be defined by Formula
$(18.1)$ for $g_{\iota }$, where ${\iota }=1, 2, ...$. The space
$L^2(V\times \Omega , \sigma ({\cal B}(V)\times {\cal F}), {\bf
l}\times \mu , Y)=: \Gamma (Y)$ is complete, since $Y$ is the
Hilbert space over $\bf F$. The sequence $\xi _{\iota }$ is
fundamental in $\Gamma (Y)$, hence converges to some $\eta $ in it.
On the other hand,
$$E(\int_V |\eta (\tau ) - \xi _{\iota }(\tau )|_Y^2{\bf l}(d\tau )) =
$$ $$ \int_V\int_{\Lambda } (\hat{{\sf m}}(d\lambda )(g(\tau , \lambda
)-g_{\iota }(\tau , \lambda )),(g(\tau , \lambda )-g_{\iota }(\tau ,
\lambda )))_X{\bf l}(d\tau )$$ for each ${\iota }$. Taking the limit
when ${\iota }$ tends to the infinity we deduce that $E(|\eta (\tau
)-\xi (\tau )|_Y^2)=0$ for $\bf l$-a.e. $\tau $ in $V$. Modifying
$\xi $ in the following manner $\hat{\xi }(\tau )=\eta (\tau )$ if $
\mu \{ \omega \in \Omega : ~ \xi (\tau )(\omega )\ne \eta (\tau
)(\omega ) \} =0$, also $\hat{\xi }(\tau )=\xi (\tau )$ if $ \mu \{
\omega \in \Omega : ~ \xi (\tau )(\omega )\ne \eta (\tau )(\omega )
\} >0$, we get a $(\sigma ({\cal B}(V)\times {\cal F}),B(Y))$
measurable random function $\hat{\xi }(\tau )$ with values in $Y$.
Two random functions $\hat{\xi }$ and $\xi $ differ on a set of
${\bf l}\times \mu $ measure null, consequently, $\hat{\xi }$ and
$\xi $ are randomly equivalent.

\par {\bf 20. Remark.} Using Lemma 19 we shall consider measurable
random vector valued functions defined by the integral like in
Formula $(18.1)$.

\par {\bf 21. Lemma.} {\it Suppose that the conditions of Remark 18 are
fulfilled. If $g\in L^2(V\times \Lambda ,\sigma ({\cal B}(V)\times
{\cal B}), {\bf l}\times \hat{{\sf m}}, X)$ and $h\in L^2(V,{\cal
B}(V), {\bf l}, {\bf F})$, then
$$(21.1)\quad \int_V \int_{\Lambda } H(d\lambda ) h(\tau )g(\tau
,\lambda ){\bf l}(d\tau )= \int_{\Lambda } H(\lambda )f(\lambda ),$$
where $f(\lambda )=\int_Vh(\tau )g(\tau ,\lambda ){\bf l}(d\tau )$.}
\par {\bf Proof.} Equality $(21.1)$ is valid for each step function
$g$. Notice that $$(21.2)\quad E(\int_{\Lambda } H(d\lambda
)f(\lambda ),\int_{\Lambda } H(d\lambda )f(\lambda ))_Y=
\int_{\Lambda }(\hat{{\sf m}}(d\lambda )f(\lambda ),f(\lambda
))_X.$$ From the Cauchy-Bunyakovskii inequality and Lemma 8 it
follows that
$$(21.3)\quad E(|\int_V\int_{\Lambda } H(d\lambda )h(\tau )g(\tau, \lambda ){\bf
l}(d\tau )|_Y^2)\le $$ $$(\int_V|h(\tau )|^2{\bf l}(d\tau ))^2 \cdot
(\int_V\int_{\Lambda } (\hat{{\sf m}}(d\lambda )g(\tau ,\lambda ),
g(\tau , \lambda ))_X{\bf l}(d\tau ))^2.$$ Taking a sequence of step
functions $g_{\iota }$ converging to $g$ in $L^2(V\times \Lambda
,\sigma ({\cal B}(V)\times {\cal B}), {\bf l}\times \hat{{\sf m}},
X)$ and using the equality $(21.2)$ and the inequality $(21.3)$ we
infer Formula $(21.1)$. \par Below applications to PDEs of
orthogonal random operator valued measures are described.

\par {\bf 22. Remark. PDEs.} Henceforth the unitary
space $\Lambda = {\bf C}^{m+3}$ is considered together with a
$\sigma $-algebra ${\cal B}={\cal B}({\Lambda })$ generated by some
semi-ring $\cal M$ of sets contained in $\Lambda $, where $m$ is the
degree of the polynomial $Q(t)$ (see Subsection 4). Take an
orthogonal random operator valued measure $H: {\cal B} \to
L^2(\Omega ,{\cal F}, \mu ; {\cal L}(X,Y)_{\kappa })$ such that
${\cal B}$ is $\hat{\sf m}$ complete, where $(\Omega , {\cal F}, \mu
)$ is a probability space, either $\tau _{\kappa }=\tau _{|\cdot |}$
or $\tau _{\kappa }=\tau _s$ (see Subsection 7). It also will be
supposed that a singleton $\{ \lambda \} $ is  in ${\cal B}$ for
each $\lambda $ in $\Lambda $. \par We consider the Sobolev space
$W_{2,m,4}([0,T]\times V_1^2, {\cal B}_{\bf l}, {\bf l}, {\bf C})$
of all complex valued functions $m$ times Sobolev in $t$ and $4$
times Sobolev in the variables $x_1,...,x_n, y_1,...,y_n$,
$x=(x_1,....,x_n)\in V_1$ and $y\in V_1$, where $V_1$ is a domain in
${\bf R}^n$, ${\bf l}$ is the Lebesgue measure restricted on
$[0,T]\times V_1^2$ from that of on ${\bf R}^{2n+1}$, ${\cal B}_{\bf
l}$ denotes an $\bf l$ complete $\sigma $-algebra on $[0,T]\times
V_1^2$ (see also Ch. III, Section 4 in \cite{mikhb}). That is
$$|f|_{W_{s,m,k}} = (\sum_{m_0=0}^m \sum_{0\le m_1+m_2\le k}
\int_0^T\int_{V_1}\int_{V_1} |\frac{\partial
^{m_0+m_1+m_2}}{\partial t^{m_0}\partial x_j^{m_1}\partial
y_k^{m_2}} f(t,x,y)|^s dt dx dy)^{1/s} <\infty $$ for each $f\in
W_{s,m,k}([0,T]\times V_1^2, {\cal B}_{\bf l}, {\bf l}, {\bf C})$,
where $1\le s<\infty $, $k\in {\bf N}$. We omit ${\cal B}_{\bf l}$
and ${\bf l}$ in order to shorten the notation. \par Let $0<T<\infty
$ and a set $V_1$ be canonically closed and compact in ${\bf R}^n$
such that to it a domain $V$ in ${\cal A}_r$ corresponds by Formula
$(1.4)$. Let also $w_0\in Int (V)$ (see $w_0$ in Formula $(5.6)$).
The Lebesgue measure on $[0,T]\times V_1^2$ induces the Lebesgue
measure on $[0,T]\times V^2$. Then Formula $(1.8)$ provides the
Sobolev space $W_{2,m,4}([0,T]\times V^2, {\bf C})$, which is taken
as a Hilbert space $X$. Then $L^2([0,T]\times V^2, {\bf C})$ is
chosen as a Hilbert space $Y$. So it can be taken a restriction from
$\cal U$ onto $V$, $~V\subset \cal U$, where $\cal U$ is a domain as
in Subsection 5.
\par Let a function $\phi (t,\lambda )$ in $t$ with a parameter
$\lambda $ be a solution of the Cauchy problem
$$(22.1)\quad Q(\frac{d}{dt}) \phi (t, \lambda ')= \lambda _1 \phi ^2(t, \lambda ')\mbox{ for }t\ge
0,$$
$$(22.2)\quad \phi (0, \lambda ')=\lambda _2, \quad \frac{d}{dt} \phi (t, \lambda ')|_{t=0}= \lambda
_3, ..., \frac{d^{m-1}}{dt^{m-1}} \phi (t, \lambda ')|_{t=0}=
\lambda _{m+1},$$ where $\lambda ' = (\lambda _1,...,\lambda
_{m+1})$, $~\lambda _1\ne 0$, $ ~\lambda _j$ is a real or a complex
constant for each $j=1,...,m+1$. In view of Theorem 1 in Subsection
3.1.5 in \cite{matvb} this Cauchy problem has a unique solution
$\phi $ belonging to $C^m([0, \infty ), {\bf C})$ in the variable
$t\in [0, \infty )$, where $C^m([0, \infty ), {\bf C})$ denotes the
space of all $m$ times continuously differentiable functions from
$[0, \infty )$ into $\bf C$. For each $0<T<\infty $ the mapping
${\bf C}^{m+1}\ni \lambda ' \mapsto \phi (\cdot, \lambda ')\in
C^m([0,T],{\bf C})$ is continuous. If $\lambda ' \in {\bf R}^{m+1}$
and $\lambda _1>0$, then $\phi $ is real valued. To simplify
notations we write $\phi (t, \lambda )$ instead of $\phi (t, \lambda
')$.
\par Let $a_1(\lambda ) = - \lambda _1$, $ ~ a_2(\lambda ) = - \alpha \lambda
_1$ and $a_3(\lambda )=\beta \lambda _1$, $~q_1(\lambda )=2\lambda
_1\lambda _{m+2}$, $~q_2(\lambda )=2\lambda _1\lambda _{m+3}$, $ ~
\lambda _{m+2}=p_1$, $ ~ \lambda _{m+3}=p_2$ and $p_2 \gamma = p_1
\varsigma $, where $\gamma $ and $\varsigma $ belong to ${\bf C}$,
$|\gamma |+|\varsigma |>0$; $p_1=0$ if $\gamma =0$, $p_2=0$ if
$\varsigma =0$, $p_1\ne 0$ if $\gamma \ne 0$, $p_2\ne 0$ if
$\varsigma \ne 0$ (see Formulas $(4.2)$ and $(4.5)$ and $(5.2)$
also). \par
\par Let $X_0=C^{m,4}([0,T]\times V^2, {\bf C})$ denote the Banach
space of $m$ times continuously differentiable functions
$f(t,z(x),z(y))$ in $t\in [0,T]$ and $(x_1,...,x_n)\in V_1$ and
$(y_1,...,y_n)\in V_1$, where $z(x)=x_1i_1+...+x_ni_n$ (see also
Formulas $(1.3)$ and $(1.4)$). The mapping ${\bf C}^{m+3}\ni \lambda
\mapsto F_{{\bf a}(\lambda )}\in C^4(V^2,{\bf C})$ is continuous
(see, for example, \cite{hormb3v,grubbmz95,mikhb}), consequently,
the mapping ${\bf C}^{m+3}\ni \lambda \mapsto K_{{\bf a}(\lambda ),
{\bf q}(\lambda )}\in C^4(V^2,{\bf C})$ also is continuous by
Formulas $(5.5)$ and $(5.6)$. Then we put $$(22.3)\quad
u(t,x,y;\omega ) = \int_{\Lambda } H(d\lambda )(\omega ) \phi (t,
\lambda ) K_{{\bf a}(\lambda ), {\bf q}(\lambda )}(x,y).$$  It will
be supposed that
\par $(22.4)$ $H(\lambda ) (\omega ) $ commutes almost $\hat{{\sf
m}} \times \mu $-everywhere with $\frac{\partial ^{m_0}}{\partial
t^{m_0}}$ and $\sigma _x^{m_1}$ and $\sigma _y^{m_1}$ and $bI$ for
every $m_0=1,...,m$ and $m_1=1,...,4$
and $b$ in ${\cal A}_{r,C}$, \\
where $\lambda \in \Lambda $, $~\omega \in \Omega $. The space $X_0$
is dense in $X$ (see the notation above). If $f$ and $g$ belong to
$X_0$, then $fg\in X_0$, where $(fg)(t,x,y)=f(t,x,y)g(t,x,y)$ for
each $(t,x,y)\in [0,T]\times V^2$. Henceforward it also will be
supposed that for each $f\in X_0$
\par $(22.5)$ $E((H(d\lambda )(\omega )f) (H(d\vartheta ) (\omega )f)) = \delta (\lambda  -\vartheta )
\xi (\lambda ) E(H(d\lambda ) (\omega )f^2)$ \\
$\hat{{\sf m}}$ almost everywhere, where $\lambda $ and $\vartheta $
are in $\Lambda $, $~\delta (\lambda  -\vartheta )$ denotes the
$\delta $ function relative to $EH(d\lambda ) (\omega )$. If $\gamma
\ne 0$ it will be taken $\xi (\lambda ) = \gamma /(2\lambda
_1\lambda _{m+2})$ for each $\lambda $ in ${\bf C}^{m+3}$ with
$\lambda _1\lambda _{m+2}\ne 0$. If $\varsigma \ne 0$ it will be
chosen $\xi (\lambda ) = \varsigma /(2\lambda _1\lambda _{m+3})$ for
each $\lambda $ in ${\bf C}^{m+3}$ with $\lambda _1\lambda _{m+3}\ne
0$ (see the notations in Subsections 4 and 5). Let
$$(22.6)\quad \rho _s(\frac{\partial ^l}{\partial t^{m_0}\partial x_j^{m_1}\partial
y_k^{m_2}}u(t,x,y;\omega )) = $$ $$E( \int_{\Lambda }
|\frac{\partial ^l}{\partial t^{m_0}\partial x_j^{m_1}\partial
y_k^{m_2}} H(d\lambda )(\omega ) \phi (t, \lambda ) K_{{\bf
a}(\lambda ), {\bf q}(\lambda )}(x,y)|^s) <\infty $$ converge
uniformly in $(t,x,y)$ on each canonically closed compact subset in
$[0,T]\times ({\cal U}\cap V)^2$ for every $m_0=0,1,...,m$, $m_1=0,
..., 4$, $m_2=0, ..., 4$ and $m_1+m_2\le 4$ if $s=1$;
$m_0=0,1,...,m$, $m_1=0,1$, $m_2=0,1$ and $m_1+m_2\le 1$ if $s=2$;
$l=m_0+m_1+m_2$, $j=1,...,n$, $k=1,...,n$, $t\in [0,T]$, $x$ and $y$
in $V\cap {\cal U}$, $0<T<\infty $. Hence the integral $(22.3)$
exists $\mu $-almost everywhere on $\Omega $ for each $t\in [0,T]$,
$x$ and $y$ in $V\cap {\cal U}$ and
$$(22.7)\quad \frac{\partial ^l}{\partial t^{m_0}\partial x_j^{m_1}\partial
y_k^{m_2}}u(t,x,y)=$$ $$E \int_{\Lambda } \frac{\partial
^l}{\partial t^{m_0}\partial x_j^{m_1}\partial y_k^{m_2}} H(d\lambda
)(\omega ) \phi (t, \lambda ) K_{{\bf a}(\lambda ), {\bf q}(\lambda
)}(x,y) ,$$ where
\par $(22.8)$ $u(t,x,y)=Eu(t,x,y;\omega )$ \\
(see Ch. V Section 1 in \cite{gihmskorb} and the Fubini Theorem in
\cite{bogachmtb,danschw}). Then From $(22.4)-(22.6)$, Lemmas 8, 19
and 21 we deduce that
$$(22.9)\quad Eu^2(t,x,y; \omega )=$$ $$\int_{\Lambda }\int_{\Lambda } E(
H(d\lambda )(\omega ) \phi (t, \lambda ) K_{{\bf a}(\lambda ), {\bf
q}(\lambda )}(x,y) H(d\eta )(\omega ) \phi (t, \eta ) K_{{\bf
a}(\eta ), {\bf q}(\eta )}(x,y))$$
$$=E\int_{\Lambda }H(d\lambda ) (\omega )\xi (\lambda )\phi ^2
(t, \lambda ) K^2_{{\bf a}(\lambda ), {\bf q}(\lambda )}(x,y).$$ On
the other hand, from $(22.5)$, $(5.5)$ and $(5.6)$ the equality
follows
$$(22.10)\quad H(d\lambda )(\omega ) \phi (t, \lambda ) K_{{\bf a}(\lambda ), {\bf q}(\lambda
)}(x,y)=H(d\lambda ) (\omega ) (\phi (t, \lambda ) F_{{\bf
a}(\lambda )}(x,y))+$$ $$ p_1\pi _1(\mbox{ }_{\sigma }\int_{w_0} ^x
(\mbox{ }_{\sigma }\int_{w_0}^y (\mbox{ }_{\sigma }\int_w^{\infty
}H(d\lambda )(\omega ) F_{{\bf a}(\lambda )}(z,v) \phi (t, \lambda )
K_{{\bf a}(\lambda ), {\bf q}(\lambda )}(w,z))dz)dv)dw)+$$
$$  p_2 \mbox{
}_{\sigma }\int_{w_0} ^y(\mbox{ }_{\sigma }\int_{w_0} ^x (\mbox{
}_{\sigma }\int_{w_0}^s (\mbox{ }_{\sigma }\int_w^{\infty
}H(d\lambda )(\omega ) F_{{\bf a}(\lambda )}(z,v) \phi (t, \lambda )
K_{{\bf a}(\lambda ), {\bf q}(\lambda )}(w,z))dz)dv)dw)ds, $$  and
$\{ Q(\frac{\partial }{\partial t}) {\bf S}_0 \} (\phi (t, \lambda )
F_{{\bf a}(\lambda )}(x,y))=0 $. \par Therefore utilizing Condition
$(22.4)$ and Formulas $(22.9)$, $(22.10)$, $(5.1)$ and
$(5.9)-(5.11)$ we deduce that $u(t,x,y; \omega )$ fulfills the PDE
$(4.4)$.
\par Therefore from Theorem 6, Lemmas 8, 19 and
21 we infer the following result.

\par {\bf 23. Theorem.} {\it If conditions $(5.8)$ and $(22.4)-(22.6)$
are fulfilled, then the random function $u(t,x,y;\omega )$ defined
by Formula $(22.3)$ is a solution of the PDE $(4.4)$ on $[0,T]\times
({\cal U}\cap V)^2$ and on the diagonal $x=y$ it satisfies the PDE
$(4.2)$.}

\par {\bf 24. Theorem}. Solutions of the Sobolev-Burgers PDE.
\par {\it Let $u(t,x)\in W_{2,m,4}([0,T]\times V_1,{\bf C})$ be a solution of the Sobolev-Burgers PDE
$(4.1)$ on $[0,T]\times V_1$, where $0<T<\infty $, $~V_1$ is a
canonically closed compact subset in ${\bf R}^n$. Then an orthogonal
random operator valued measure $H$ exists such that
$u(t,x)=u(t,z(x),z(x))$, where $u(t,x,y)$ is given by Formulas
$(22.3)$ and $(22.8)$, $z(x)=x_1i_1+...+x_ni_n$.}
\par {\bf Proof.} {\bf I.} If in addition to conditions of Theorem
23
\par $(24.1)$ $E(u^2(t,x,y; \omega ))=u^2(t,x,y)$, \\ where $u(t,x,y)$ is provided
by Formulas $(22.3)$ and $(22.8)$, then $u(t,x)=u(t,z(x),z(x))$ is a
solution of the PDE $(4.1)$ (see also Formulas $(1.3)$ and $(1.4)$).
\par From Subsections 4, 5 and 22 it follows that a domain \par
$(24.2)$ $\Upsilon := \{ \lambda \in {\bf C}^{m+3}: \exists ~
v(t,x)=\phi (t, \lambda ) K_{{\bf a}(\lambda ), {\bf q}(\lambda
)}(z(x),z(x))$\par $ \mbox{ satisfying the PDE (4.1) on }
[0,T]\times  V_1 \} $\\
is a Borel set in ${\bf C}^{m+3}$. We mean here the correspondence
between PDEs in variables belonging to domains in ${\bf R}^n$ and in
${\cal A}_r$ according to Subsections 1 and 4 and simplifying the
notation.
\par {\bf II.} It is sufficient to provide $H(\lambda ) (\omega )$
on $\Upsilon \times \Omega $ and extend it by zero on $({\bf
C}^{m+3}\setminus \Upsilon )\times \Omega $. We consider a family
$\cal Z$ of countable disjoint partitions $\Phi = \{ G_1, G_2,... \}
$ of $\Lambda = {\bf C}^{m+3}$ possessing the following properties
\par $(24.3)$ $\bigcup_{\iota =1}^{\infty }G_{\iota }=\Lambda $,  $ ~ G_{\iota }\cap G_j=\emptyset
$ for each ${\iota }\ne j$, $~\sup_{\iota } (diam (G_{\iota }))=:
d(\Phi )<\infty $, $ ~ cl (G_{\iota }) = cl (Int(G_{\iota }))$, $ ~
G_{\iota }$ belongs to the Borel $\sigma $-algebra $B(\Lambda )$ of
$\Lambda $, where $Int (G_{\iota })$ denotes the interior of
$G_{\iota }$, whilst $cl (G_{\iota })$ denotes the closure of
$G_{\iota }$, $ ~ diam (G_{\iota })= \sup_{a, b\in G_{\iota }}
|a-b|$, $ ~ |a|=\sqrt{|a_1|^2+...|a_{m+3}|^2}$ for each
$a=(a_1,...,a_{m+3})\in {\bf C}^{m+3}$. The family $\cal Z$ is
partially ordered $\Phi _1=\{ G_{1,1}, G_{1,2},... \} <\Phi _2= \{
G_{2,1}, G_{2,2},... \} $ if and only if for each $G_{2,k}$ a
natural number $j(k)$ exists such that $G_{2,k}\subset G_{1,j(k)}$,
that is a partition $\Phi _2$ is finer than $\Phi _1$.
\par Evidently, the linear span $span _{\bf C} \{ \phi (\cdot ,
\lambda ): ~ \lambda \in \Lambda \} $, where $\phi (\cdot , \lambda
)$ is a solution of the Cauchy problem $(22.1)$ and $(22.2)$ is
dense in $W_{2,m}^Q:=\{ g\in W_{2,m}([0,T],{\bf C}): ~ g\notin {\cal
N} (Q(\frac{d}{dt})) \} $, where ${\cal N}(Q(\frac{d}{dt}))$ denotes
the null-space of the differential operator $Q(\frac{d}{dt})$.
Indeed, $\phi (t, \lambda )$ is in $C^m([0,T],{\bf C})$ for each
$\lambda \in \Lambda $. A set $\{ \lambda _1f^2: ~ f\in
C^m([0,T],{\bf R}), ~ \lambda _1>0 \} $ has a non void interior in
$C^m([0,T],{\bf R})$. On the other hand, $C^{2m}([0,T],{\bf R})$ and
$Q(\frac{d}{dt})C^{2m}([0,T],{\bf R})$ are dense in $C^m([0,T],{\bf
R})$. Therefore, $\Psi := \hat{\pi }_2 \{ (f, \lambda
_1f^2)=(g,Q(\frac{d}{dt})g): ~ f\in C^m([0,T],{\bf R}), ~ \lambda
_1>0, ~g\in C^{2m}([0,T],{\bf R}) \} $ is dense in some open subset
in $\{ g\in C^m([0,T],{\bf R}): g\notin {\cal N}(Q(\frac{d}{dt})) \}
$, where $\hat{\pi }_2$ denotes the linear projection $\hat{\pi
}_2(f,h)=h$ for all $f$ and $h$ in $C^m([0,T],{\bf R})$. The space
$C^m([0,T],{\bf C})$ is dense in $W_{2,m}([0,T],{\bf C})$, hence the
$\bf C$ linear span of $\Psi $ is dense in $W_{2,m}^Q$.
\par The Sobolev spaces $W_{2,m,4}([0,T]\times V^2,{\bf C})=X$,
$W_{2,m}([0,T],{\bf C})$ and $W_{2,4}(V^2,{\bf C})$ are separable.
Suppose that $u(t,x,y)$ is in $W_{2,m,4}([0,T]\times V^2,{\bf C})$
and $u(t,x,x)$ is a solution of the PDE $(4.1)$. If $f\in X$, then
$Q(\frac{\partial }{\partial t}){\bf S}_0f\in Y$, where
$Y=L^2([0,T]\times V^2, {\bf C})$. On the other hand, the tensor
product $L^2([0,T],{\bf C})\otimes L^2(V^2,{\bf C})$ is dense in
$Y$. For each $\epsilon
>0$ there are non null functions $\phi _j(t)=\phi (t,\lambda (j))$ and
$f_j(x,y)\in W_{2,4}(V^2,{\bf C})$ such that $\lambda (j)\in \Lambda
$ for each $j=1,...,k$ and $|u-\sum_{j=1}^k\phi _jf_j
|_{W_{2,m,4}}<\epsilon $ and $\lambda (j)\ne \lambda ({\iota })$ for
each ${\iota }\ne j$, also functions $\{ v_j:=\phi _jf_j: j=1,...,k
\} $ are linearly independent, where $k=k(\epsilon )\in \bf N$. Thus
$v_j\in X\subset Y$ for each $j$.
\par Obviously $X$ and $Y$ considered as the $\bf R$ linear spaces
are isomorphic with $X_R\oplus X_R{\bf i}$ and $Y_R\oplus Y_R{\bf
i}$ respectively, where $X_R=W_{2,m,4}([0,T]\times V^2,{\bf R})$ and
$Y_R=L^2([0,T]\times V^2,{\bf R})$. Let $\hat{\pi }_{v_1,...,v_k}:
Y\to span_{\bf C}(v_1,...,v_k)$ be a $\bf C$-linear projection
operator, where \par $span_{\bf C}(v_1,...,v_k)= \{ v\in Y: ~
v=c_1v_1+...+c_kv_k: ~ c_1\in {\bf C},...,c_k\in {\bf C} \} $ \\
denotes the linear span of vectors $v_1,...,v_k$ over $\bf C$. Since
$v_1,...,v_k$ are linearly independent, then $span_{\bf
C}(v_1,...,v_k)$ is isomorphic with ${\bf C}^k$. If it is considered
as the $\bf R$ linear space, then it is isomorphic with ${\bf
R}^{2k}$. Then a partition $\Phi =\Phi (v_1,...,v_k)$ with $\Phi \in
\cal Z$ exists so that \par $(24.4)$ for each $1\le j\ne l\le k$
there are unique ${\iota }(j)\ne {\iota }(l)$ with $\lambda (j)\in
G_{{\iota }(j)}$ and $\lambda (l)\in G_{{\iota }(l)}$. \par Let
${\cal M}={\cal M}(v_1,...,v_k)$ be a semi-ring of subsets in
$\Lambda $ generated by a disjoint partition $\Phi = \{ G_1, G_2,
... \} $ satisfying Conditions $(24.3)$ and $(24.4)$.
\par  {\bf III.} Take a probability space
$(\Omega _{v_1,..,v_k}, {\cal G}_{v_1,..,v_k}, P_{v_1,..,v_k})$ with
\par $\Omega _{v_1,..,v_k}=span_{\bf C}(v_1,...,v_k)$, where
$P_{v_1,..,v_k}: {\cal G}_{v_1,..,v_k}\to [0,1]$ is a probability,
${\cal G}_{v_1,..,v_k}$ is a completion of the Borel $\sigma
$-algebra $B(\Omega _{v_1,..,v_k})$ relative to $P_{v_1,..,v_k}$. We
put $v(\lambda )=v_j$ for each $\lambda \in G_{\iota (j)}$ Choose an
orthogonal random operator valued measure $H_{v_1,..,v_k}(d\lambda
)$ and a probability measure $P_{v_1,..,v_k}$ such that to satisfy
$(22.4)$, $(22.5)$ and $(24.1)$ for $H_{v_1,..,v_k}(d\lambda )$
restricted on $span_{\bf C}(v_1,...,v_k)$ and for
$$u_{v_1,...,v_k}=\int_{\Lambda }H_{v_1,..,v_k}(d\lambda )v(\lambda
)$$ instead of $u$. This $P_{v_1,..,v_k}$ can be taken as a Gaussian
measure corresponding to a random vector $\upsilon _{v_1,...,v_k}$
with a mean value $\theta _{v_1,...,v_k}$ and a correlation operator
${\cal C}_{v_1,...,v_k}$.
\par The Hilbert spaces $X$ and $Y$ are separable,
consequently, they are isomorphic. If $\hat{J}: Y\to X$ is an
isomorphism, then $\hat{J}H: X\to X$. By virtue of the Fubini
theorem Conditions $(22.5)$ and $(24.1)$ mean that $\hat{H}^2=\Xi
\hat{H}$, where $\hat{H}=\hat{J}EH$. Here an operator $\Xi $ is the
multiplication operator on a function $\xi (\lambda )$, $ ~ \Xi
y(\lambda ) = \xi (\lambda )y(\lambda )$ for each $\lambda \in
\Lambda $ and $y\in L^2(\Lambda , {\cal B}, \hat{{\sf m}}, X)$. In
order to describe an orthogonal random operator valued measure
$\hat{J}H_{v_1,..,v_k}(M)$ for each $M\in {\cal M}(v_1,...,v_k)$ it
is sufficient to provide $\hat{J}H_{v_1,..,v_k}(G_j)$ for each $G_j$
of the disjoint partition $\Phi = \{ G_1, G_2, ... \} $. For each
$j=1, 2,...$ an orthogonal random operator valued measure
$\hat{J}H_{v_1,..,v_k}(G_j)$ can be realized as a random normal
operator and considering its representation with the help of an
integral by a random projection operator valued measure on its
spectrum (see Section 5.2 in \cite{kadring}, also
\cite{danschw,gihmskorb}).
\par Then we take a monotone decreasing sequence $\epsilon
_1>\epsilon _2>...>0$ and choose $v_1, v_2,...$ such that
$$(24.5)\quad |u-\sum_{j=1}^{k(\iota )}v_j|_{W_{2,m,4}}<\epsilon _{\iota }$$ for
each $\iota =1, 2,...$, where $k(\iota )<k(j)$ for each $\iota <j$.
Then for each $j$ a countable partition $\Phi _j$ of $\Lambda $ can
be chosen satisfying Conditions $(24.3)$ and $(24.4)$ and such that
$\Phi _{\iota }<\Phi _j$ for each $\iota <j$ and $\lim_{j\to \infty
}d(\Phi _j)=0$.
\par Therefore it can be put $\Omega = cl_Y ( span_{\bf C}(v_j: j \in
{\bf N} )) $, where the closure of a linear span is taken in $Y$.
There are natural projections $\hat{\pi }_k: \Omega \to \Omega
_{v_1,..,v_k}$ and $\hat{\pi }^l_k: \Omega _{v_1,..,v_l}\to \Omega
_{v_1,..,v_k}$ for each $l>k$.
\par A probability Gaussian measure $P_{v_1,..,v_k}$ on ${\cal
G}_{v_1,..,v_k}$ induces a cylindrical measure $\mu _{v_1,...,v_k}$
on a cylindrical $\sigma $-algebra ${\cal F}(v_1,...,v_k)$ of all
cylindrical sets $(\hat{\pi }_{v_1,...,v_k})^{-1}(N)$ in $\Omega $
with $N$ in ${\cal G}_{v_1,..,v_k}$. This family of measures $\mu
_{v_1,...,v_k}$ can be chosen consistent with projections so that it
induces a bounded cylindrical distribution on an algebra $\bigcup_k
{\cal F}(v_1,...,v_k)$.  Mean values $\theta _{v_1,...,v_k}$ and
correlation operators ${\cal C}_{v_1,...,v_k}$ can be chosen so that
$\sup_k |\theta _{v_1,...,v_k}|_X<\infty $ and $\sup_k |{\cal
C}_{v_1,...,v_k}| <\infty $, where $|{\cal C}_{v_1,...,v_k}|$
denotes a norm of ${\cal C}_{v_1,...,v_k}$ for $span_{\bf
C}(v_1,...,v_k)$ embedded into $(X, |\cdot |)$. The natural
embedding operator ${\cal J}_0: X\hookrightarrow Y$ is nuclear (i.e.
of trace class). Therefore, the bounded consistent family of
measures induces a $\sigma $-additive Gaussian measure $\mu $ on
$(\Omega , {\cal F})$, where ${\cal F}$ is a completion of
$\bigcup_k {\cal F}(v_1,...,v_k)$ (see Section II.2 in
\cite{dalfb}).
\par Then the family of orthogonal random operator valued measures
$H_{v_1,..,v_k}(d\lambda )$ on subspaces $span_{\bf C}
(v_1,...,v_k)$ and subalgebras ${\cal B}_k$ generated by ${\cal
M}(v_1,...,v_k)$ induces an orthogonal random operator valued
measure $H(d\lambda )$ on $X(v_j: j\in {\bf N}):=cl_X ( span_{\bf
C}(v_j: j \in {\bf N} )) $, since $\mu ((\hat{\pi
_k})^{-1}(N))=\mu_{v_1,...,v_k}(N)$ for each $N\in {\cal
G}_{v_1,..,v_k}$. It has an extension by the identity operator $I$
on the orthogonal complement $X\ominus X(v_j: j\in {\bf N})$.
Naturally ${\cal B}$ is the completion of $\bigcup_k {\cal B}_k$,
where ${\cal B}_k$ denotes a completion of ${\cal M}(v_1,...,v_k)$,
$k=k(\iota )$, $ ~ \iota =1, 2,...$. This provides properties
$(22.4)$, $(22.5)$ and $(24.1)$ for $H(d\lambda )$ on $X$. Taking a
limit of a fundamental sequence of step functions in $L^2(\Lambda ,
{\cal B}, \hat{\sf m}, {\bf C})$ gives a function $v\in L^2(\Lambda
, {\cal B}, \hat{\sf m}, {\bf C})$ (see Subsection 10), where
$v(\lambda )=\phi (\cdot ,\lambda )f(\cdot , \lambda )$ with
$f(\cdot ,\lambda )\in W_{2,4}(V^2,{\bf C})$ for each $\lambda $,
also $v(\lambda (j))=v_j$ for each $j$. By virtue of Lemma 21 we get
$$u(t,x,y;\omega )= \int_{\Lambda }H(d\lambda )(\omega )v(\lambda
).$$
\par Using Formulas $(5.5)$ and $(5.6)$ define a function
$g(x,y,\lambda )$ related with $f(x,y,\lambda )$ similarly to the
pair of $F_{{\bf a}(\lambda )}$ and $K_{{\bf a}(\lambda ), {\bf
q}(\lambda )}$, where $\lambda $ is the parameter in $\Lambda $ so
that $\hat{\sf m}(G)\ne 0$ for each $j$ and $G$ in $\Phi _j$ with
$\lambda \in G$. Therefore, utilizing these properties together with
the orthogonality of $H$ we infer for a function $h(t,x,y,\lambda )=
\phi (t,\lambda )g(x,y,\lambda )$ that $$E(|\int_{\Lambda
}H(d\lambda )Q(\frac{\partial }{\partial t}){\bf S}_0h(t,x,y,
\lambda )|_{x=y}^2)=E(|Q(\frac{\partial }{\partial t}){\bf
S}_0\int_{\Lambda }H(d\lambda )h(t,x,y,\lambda )|_{x=y}^2)=0$$ $\bf
l$ almost everywhere on $[0,T]\times V$, where $a=-\alpha $ and
$b=\beta $ (see Formulas $(4.1)$ and $(4.5)$). This implies that
$v(\lambda )(t,x,x)=\phi (t ,\lambda )f(x,x, \lambda )$ satisfies
the PDE $(4.1)$ for $\hat{\sf m}$ almost all $\lambda $ in $\Lambda
$. Thus $f(\cdot , \lambda )= K_{{\bf a}(\lambda ), {\bf q}(\lambda
)}$ for $\hat{\sf m}$ almost all $\lambda $ in $\Lambda $.

\par {\bf 25. Remark.} The Laplacian $\Delta _x$ is invariant under each orthogonal
transformation $T\in O(n)$, where $O(n)$ notates the orthogonal
group of the Euclidean space ${\bf R}^n$. If there is a vector
Sobolev-Burgers PDE with $u\in {\bf R}^k$ and the scalar product
$(u,u)$ instead of $u^2$, where $k\le 4$, it is possible to consider
its noncommutative analog for $u$ with values in the quaternion skew
field ${\bf H}={\cal A}_2$. The following generalized PDE
$$(25.1)\quad E \{ Q(\frac{\partial }{\partial t}) {\bf S}_0 u(t,x,y;\omega )
+ (\sigma _x+\sigma _y)(u^2(t,x,y;\omega ) q_1 + u^2(t,x,y;\omega )
q_2 \} |_{x=y}=0$$ can be solved analogously to Subsection 4-6, 22,
23, where $q_1$ and $q_2$ are in the quaternion skew field ${\cal
A}_2={\bf H}$ and $|q_1|+|q_2|>0$, $x$ and $y$ are in $V$, $2\le
r\le 3$ and $2^r>n$. For this we put also
\par $(25.2) \quad K(x,y)=F(x,y)+{\sf A}K(x,y)$ with
$$(25.3)\quad {\sf A}K(x,y)= (\mbox{
}_{\sigma }\int_{w_0} ^x (\mbox{ }_{\sigma }\int_{w_0}^y (\mbox{
}_{\sigma }\int_w^{\infty } F(z,v)K(w,z)dz)dv)dw)p_1$$ $$+ (\mbox{
}_{\sigma }\int_{w_0} ^y(\mbox{ }_{\sigma }\int_{w_0} ^x (\mbox{
}_{\sigma }\int_{w_0}^s (\mbox{ }_{\sigma }\int_w^{\infty }
F(z,v)K(w,z)dz)dv)dw)ds) p_2,$$ where $p_j\ne 0$ if and only if
$q_j\ne 0$, $j\in \{ 1, 2 \} $, $ ~ K=K_{{\bf a}, {\bf q}}$; $ ~
F=F_{\bf a}$; $ ~ w_0$ is marked point in ${\cal U}$, ${\cal
U}\subset {\cal A}_3$, $x$ and $y$ are in $V$ with $V\subset {\cal
U}$. Formulas for solutions are similar, since the octonion algebra
${\cal A}_3={\bf O}$ is alternative, while $K$ and $F$ have values
in the quaternion skew field $\bf H$.

\par {\bf 26. Conclusion.} In this paper the Sobolev-Burgers
PDEs were integrated using noncommutative line integral over
Cayley-Dickson algebras and orthogonal random operator valued
measures. The nonlinear problem was reduced to the linear PDEs for
$F$ and the $\bf C$-linear integral equation relating $F$ and $K$.
\par It is worth to mention that algorithms for numerical solutions of
integral equations converge better than for PDEs. It extends
previous approaches based on real and complex numbers, because each
scalar or vector PDE over them can be reformulated over octonions
and Cayley-Dickson algebras and new types  of PDEs can be
encompassed (see Section 2). \par The obtained results can be used
for further investigations of PDEs and properties of their
solutions. For example, generalized PDEs including terms such as
${\Delta }^p$ or ${\nabla }^p$, $div (|\nabla u|^p\nabla u)+\lambda
|u|^pu$ for $p>0$ can be investigated, for dynamical nonlinear
processes, air target range radar measurements
\cite{hormbl,kichenb96,kuzgadznlprtj,svalkorplb,zabokoportj16},
which have technical applications and in the sciences.

\end{document}